# Performance Comparison of Optimal Fractional Order Hybrid Fuzzy PID Controllers for Handling Oscillatory Fractional Order Processes with Dead Time

Saptarshi Das[a,b] [*], Indranil Pan[a], and Shantanu Das[c]

a) Department of Power Engineering, Jadavpur University, Salt Lake Campus, LB-8, Sector 3, Kolkata-700098, India.
b) Communications, Signal Processing and Control Group, School of Electronics and Computer Science, University of Southampton, Southampton SO17 1BJ, United Kingdom.
c) Reactor Control Division, Bhabha Atomic Research Centre, Mumbai-400085, India.

**Authors' Emails:**
saptarshi@pe.jusl.ac.in, s.das@soton.ac.uk (S. Das*)
indranil.jj@student.iitd.ac.in, indranil@pe.jusl.ac.in (I. Pan)
shantanu@magnum.barc.gov.in (Sh. Das)

**Abstract:**
Fuzzy logic based PID controllers have been studied in this paper, considering several combinations of hybrid controllers by grouping the proportional, integral and derivative actions with fuzzy inferencing in different forms. Fractional order (FO) rate of error signal and FO integral of control signal have been used in the design of a family of decomposed hybrid FO fuzzy PID controllers. The input and output scaling factors (SF) along with the integro-differential operators are tuned with real coded genetic algorithm (GA) to produce optimum closed loop performance by simultaneous consideration of the control loop error index and the control signal. Three different classes of fractional order oscillatory processes with various levels of relative dominance between time constant and time delay have been used to test the comparative merits of the proposed family of hybrid fractional order fuzzy PID controllers. Performance comparison of the different FO fuzzy PID controller structures has been done in terms of optimal set-point tracking, load disturbance rejection and minimal variation of manipulated variable or smaller actuator requirement etc. In addition, multi-objective Non-dominated Sorting Genetic Algorithm (NSGA-II) has been used to study the Pareto optimal trade-offs between the set point tracking and control signal, and the set point tracking and load disturbance performance for each of the controller structure to handle the three different types of processes.

## 1. Introduction:

Fractional order modeling and fractional order controllers have got recent popularity in the control engineering community due to its extra flexibility to understand, represent and control dynamical systems [1]-[3]. FO controllers have additional advantages of handling delicate and complicated process dynamics at the cost of increased computational cost. Also, hardware implementation of such controllers is difficult due to the infinite dimensional nature of FO differ-integrators, which is the building block of FO controllers. Recent research thrust has been focused on the implementation of FO controllers using fractance (electro-chemical), analog and digital electronic circuit realization [1], [2].

Fuzzy logic controllers (FLC) have been traditionally used for efficient control of nonlinear, time varying and vague systems with little knowledge of the process to be controlled. Conventional FLCs work with the loop error and its integer order differ-integral [4] and have been proved to be an effective means over conventional PID controllers. The FLC based PID controllers have certain advantages as it combines the potential of both FLC and conventional PID controller. Tuning of such fuzzy PID controllers can be done to meet the design specifications. The tuning parameters of a fuzzy PID controller can be the input-output scaling factors, rule base, shape and degree of overlap of the membership functions (MF) etc. It has been suggested in Woo *et al.* [5] that the output scaling factors can be considered as the effective gains of the FLC based PID controllers. Also, these SFs have the highest impact on the control performance over the other FLC parameters. So, tuning of the input-output SFs based on fixed MFs and rule-base can be a logical approach as shown in [5]-[7]. Moreover, similar to the different decomposed structure of conventional PID controller [8]-[9], fuzzy PID controllers have been categorized in this paper to produce enhanced closed loop performance. Golob [8] studied various decomposed and hybrid fuzzy PID structures and their relative merits in closed loop control design which has been extended in the present work with its fractional order enhancements.

Application of fractional calculus based methods to enhance the performance of conventional fuzzy logic based systems has been a recent focus in the contemporary research community. Efe [10] used fractional order integration while designing an Adaptive Neuro-Fuzzy Inference System (ANFIS) based sliding mode controller. Delavari *et al.* [11] proposed a fuzzy fractional sliding mode controller and tuned its parameters with GA. Tian, Huang and Zhang [12] has shown that fuzzy enhancement of FOPID controllers is more advantageous over that with simple PID controllers. In this paper, the concept of fuzzy enhancement of fractional order controllers has been extended by using various hybrid combinations with their optimal time-domain tuning. The motivation of the present work is to use the advantage of optimal $PI^\lambda D^\mu$ controllers while also enjoying the benefits of FLC in various combinations with the integro-differential actions for the control of complicated processes. The fuzzy logic controller in a closed loop control system can be visualized like a static non-linearity between its inputs and outputs, which can be tuned easily to match the desired performance of the control system in a more heuristic manner without delving into the exact mathematical description of the modeled nonlinearity. The proposed family of hybrid FO fuzzy PID controllers works with fractional order rate and integration of the error signal like the fractional order $PI^\lambda D^\mu$ controller, proposed by Podlubny [13].

The rationale of incorporating fractional order differ-integral actions before and after fuzzy inferencing, as studied in the present paper, needs extra illustration. If it be assumed that a human operator replaces the automatic controller in the closed loop feedback system, the operator would rely on his intuition, experience and practice to formulate a control strategy and he would not do the differentiation and integration in a mathematical sense. This can be viewed as the rationale behind the coupling between the FO differ-integration with the fuzzy inferencing, where rate of change and history of the error signal does not resemble the integer order calculus or conventional differentiation and integration in a pure mathematical sense. Rather it reflects an operator's experience which gives extra freedom for tuning of control loops. Thus, the control signal generated as a result of his actions may be approximated by appropriate mathematical operations which have the required compensation characteristics. The rationale behind incorporating fractional order operators in the conventional hybrid fuzzy PID input and output can be visualized like a heuristic reasoning for an observation of a particular rate of change in error (not in mathematical sense) by a human operator and the corresponding actions he takes over time which is not static in nature since the fractional differ-integration involves the past history of the integrand and as if the integrand is continuously changing over time [2]. Since, human brain does not observe the rate of change of a variable and its time evolution as classical integer order numerical differentiation and integration, the fractional order of differ-integration perhaps puts some extra flexibility to map information in a more easily decipherable form. Considering these flexibilities incorporated in the fuzzy inference input and output, the present study extends the idea with different hybrid structures of the FLC based FOPID controller and their comparative merits in closed loop control with fixed MF type and rule base and fuzzy inferencing.

The rest of this paper is organized as follows. Section 2 discusses the different hybrid FO fuzzy PID controller structures. Section 3 presents the methodology of single and multi-objective evolutionary optimization based time domain tuning of the optimal controller parameters i.e. the input-output SFs and the integro-differetial orders. Three different class of oscillatory fractional order processes have been tested with these controller structures in section 4 and the closed loop performances are also compared. The paper ends with the conclusion as section 5, followed by the references.

## 2. Structures of the family of fractional order hybrid fuzzy PID controllers

In a simple PID control scheme, a process when excited with an external input $r(t)$ to produce a response $y(t)$, the loop error $e(t) = r(t) - y(t)$ is minimized with simple proportional, derivative and integral actions of the error. Podlubny [13] proposed the concept of fractional order differ-integral actions on the error signal to tune $PI^\lambda D^\mu$ controllers that gives higher degrees of freedom in controller design. The present study extends the concept of FO differ-integral actions in various FLC based FOPID structures. These FLC based FOPID structures, have been inherited from their integer order counterparts as studied by Golob [8].

Though three-input FLC has been studied in several literatures [4], [8] but it finds lesser scope for practical implementation as the FLC mimics the knowledge base of a human operator who intuitively manipulates the control law by observing the error and its rate of change. For a human operator the observation of error and its rate of change are more practical than the observation of errors integral or its history. Thus, in the present study integral of error has not been considered as an input to the fuzzy inference system because it is hard to visualize and complicates the fuzzy inferencing, fuzzification-defuzzification systems and also development of an effective knowledge (rule) similar to that of a human operator. Golob [8] developed different rule base for the proportional, integral and derivative actions which may be computationally intensive for practical implementation. We have retained the same rule base for the two-input fuzzy inference engine in various combinations so as to provide the controller structure more flexibility for optimal time domain tuning. Different FO fuzzy PID controller structures are now proposed in the following subsections.

### 2.1. Fractional order fuzzy PID controller

The conventional fuzzy PID controller has been extended here in FO domain to provide the controller some extra flexibility as shown in Fig. 1. With integer order differ-integration, this typical structure has been extensively studied in [14]-[20] with various applications like liquid level control systems [19], random delay handling in networked control systems [20] etc. In Fig. 1, $K_e$ and $K_d$ are the input SFs and $K_{PI}$ and $K_{PD}$ are the output SFs. It is shown by Engin *et al.* [19] that the ratio of output SFs balances the impact of the integro-differential actions of the controller. Here, the integer order rate of error in the conventional integer order FLC input has been replaced by its FO counterpart ($\mu$). Also the FLC output is fractionally integrated with a flexible order ($\lambda$) which can be tuned to meet design specifications.

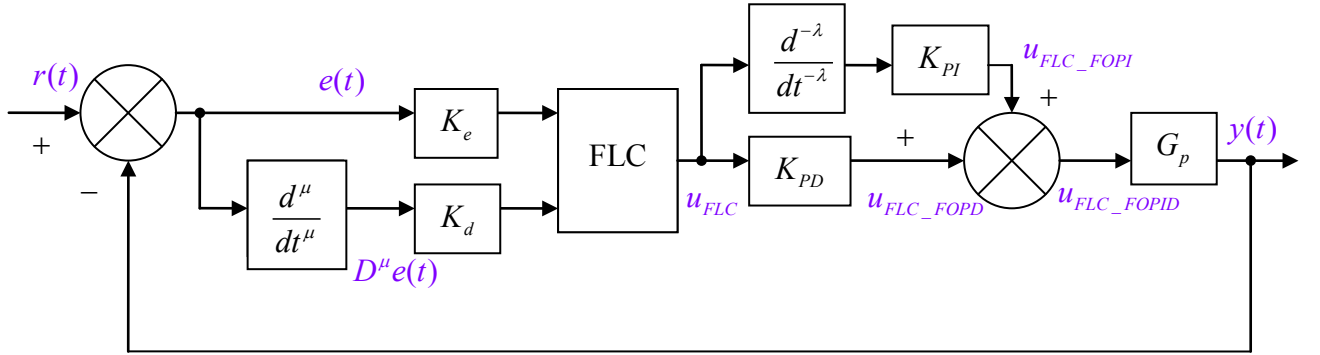

Fig.1. Schematic diagram of FO fuzzy PID controller.

This FO fuzzy PID controller can be viewed like a parallel combination of a fuzzy FOPD and fuzzy FOPI controller. The control law of the above scheme is given as:

$$u_{FLC\_FOPID}(t) = u_{FLC\_FOPI}(t) + u_{FLC\_FOPD}(t)$$
$$= K_{PI} \cdot \frac{d^{-\lambda} u_{FLC}(t)}{dt^{-\lambda}} + K_{PD} u_{FLC}(t) \qquad (1)$$

In [14] it has been shown that for this typical integer order fuzzy PID controller with product-sum inferencing, center of gravity defuzzification method and triangular membershiop function, the relation between the input and output variables can be expressed as

$$u = A + Be + C\dot{e} \qquad (2)$$

where, parameters $\{A, B, C\}$ are devised from the error and control signal and detailed in [14]. Therefore, it is expected that for the controller in Fig. 1, the FLC output (2) will be a function of fractional rate of error instead of conventional integer order derivative of error signal. Now, using the well-known identity of fractional calculus or successive derivative of higher powers of a variable $\frac{d^n}{dt^n} t^m = \frac{\Gamma(m+1)}{\Gamma(m-n+1)} t^{m-n}$, equation (1) can be expressed as

$$u_{FLC\_FOPID}(t) = K_{PI} \cdot \frac{d^{-\lambda}}{dt^{-\lambda}}\left(A + BK_e e + CK_d \frac{d^\mu e}{dt^\mu}\right) + K_{PD}\left(A + BK_e e + CK_d \frac{d^\mu e}{dt^\mu}\right)$$

$$= \left[K_{PD}A + K_{PI}A\frac{t^\lambda}{\Gamma(\lambda+1)}\right] + [K_{PD}BK_e]e + [K_{PD}CK_d]\frac{d^\mu e}{dt^\mu} \quad (3)$$

$$+ [K_{PI}BK_e]\frac{d^{-\lambda}e}{dt^{-\lambda}} + [K_{PI}CK_d]\frac{d^{\mu-\lambda}e}{dt^{\mu-\lambda}}$$

Thus drawing an analogy with the PID controller, the first term $\left[K_{PD}A + K_{PI}A\frac{t^\lambda}{\Gamma(\lambda+1)}\right]$ represents a time dependent gain due to the presence of time in it. The term $[K_{PD}BK_e]$ represents the proportional gain, $[K_{PD}CK_d]$ represents the fractional order derivative gain, $[K_{PI}BK_e]$ represents the fractional order integral gain and $[K_{PI}CK_d]$ represents an additional FO integro-differential gain. The last term can represent either a fractional derivative or a fractional integral action depending on which value between $\lambda, \mu$ is greater.

### *2.2. Fractional order fuzzy PI+PD controller*

The same rule-base can still be used for both of the FLC based FOPI and FOPD controller shown in Fig. 2, with the provision of choosing the input-output SFs independently which gives lesser complicacy in the knowledge base to implement the controller as suggested in Golob [8]. The structure shown in Fig. 2 is still a parallel combination of fuzzy FOPI and fuzzy FOPD controller with different input SFs as $\{K_{e_1}, K_{d_1}\}$ and $\{K_{e_2}, K_{d_2}\}$. This two stage structure has been studied in [21]-[22] in integer order domain. The two distinct parts of the FO fuzzy PI+PD controller uses same fractional order rate for fuzzy inferencing but after getting multiplied with different optimally tuned input SFs, which change the universe of discourse of FLC-1 and FLC-2 in a different way to give the structure more flexibility.

The control law for the FO fuzzy PI+PD controller is given by:

$$u_{FLC\_FOPID}(t) = u_{FLC\_FOPI}(t) + u_{FLC\_FOPD}(t)$$
$$= K_{PI} \cdot \frac{d^{-\lambda}u_{FLC-1}(t)}{dt^{-\lambda}} + K_{PD}u_{FLC-2}(t) \quad (4)$$

The rule base for the fuzzy FOPI and the fuzzy FOPD part is as follows:

$R_r^{PI}$ : IF $e(t)$ is $E^r$ AND $D^\mu e(t)$ is $\Delta E^r$, THEN $u_{FLC-1}$ is $\Delta U_{PI}, \forall\, r \in [1, n]$

$R_r^{PD}$ : IF $e(t)$ is $E^r$ AND $D^\mu e(t)$ is $\Delta E^r$, THEN $u_{FLC-1}$ is $U_{PD}, \forall\, r \in [1, n]$

where, $r$ represents each of the $n$ different rules of the rule base. Thus the only difference between the two parts is the control action $\Delta U_{PI}$ and $U_{PD}$.

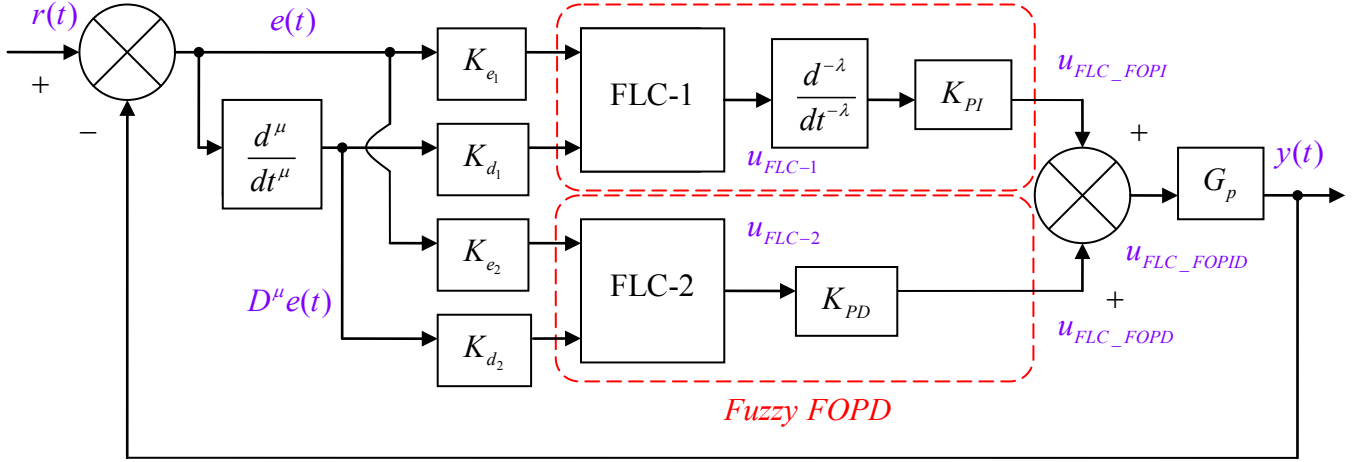

Fig. 2. Schematic diagram of FO fuzzy PI+PD controller.

### 2.3. Fractional order fuzzy P+ID controller

Li [23] proposed the concept of fuzzy P+ID controller where the fuzzy module added in conjunction with the proportional gain modulates the peak overshoot and rise time efficiently. The conventional integral and derivative actions preserve the zero-offset and enhanced stability with flatness of the output signal respectively. In this structure the derivative action is implemented on the process output like [23]-[24] and not in the conventional way. Thus derivative action is smooth in this case which prevents the derivative kick for sudden step change in set-point. This is especially needed in process control applications as the controller senses a sudden jump in error rate and to suppress it the corresponding derivative action becomes very large. The fractional derivative action with gain $K_{d_2}$ and order $\mu_2$ in the feedback path and the fractional integral action in forward path with gain $K_i$ and order $\lambda$ modulate the rate of change in process variable and time evolution of error signal respectively in a more delicate manner while producing enhanced closed loop performance in terms of handling derivative-kick for sudden step change in set-point.

This FLC structure preserves the basic simplicity of the conventional PID controller with an extra fuzzy module based proportional action which is easier to implement in real life hardware. It is also shown by Li [23] that if a PID controller gives stable response for a specific design, its fuzzy P counterpart also guarantees stability, though FLC introduces an extra nonlinearity in the design.

The control law in this case is a combination of fuzzy P and conventional FO integral-derivative controller is given by:

$$u_{FLC\_FOPID}(t) = u_{FLC\_P}(t) + u_{FOI}(t) - u_{FOD}(t)$$
$$= K_p u_{FLC}(t) + K_i \cdot \frac{d^{-\lambda} e(t)}{dt^{-\lambda}} - K_{d_2} \cdot \frac{d^{\mu_2} y(t)}{dt^{\mu_2}} \quad (5)$$

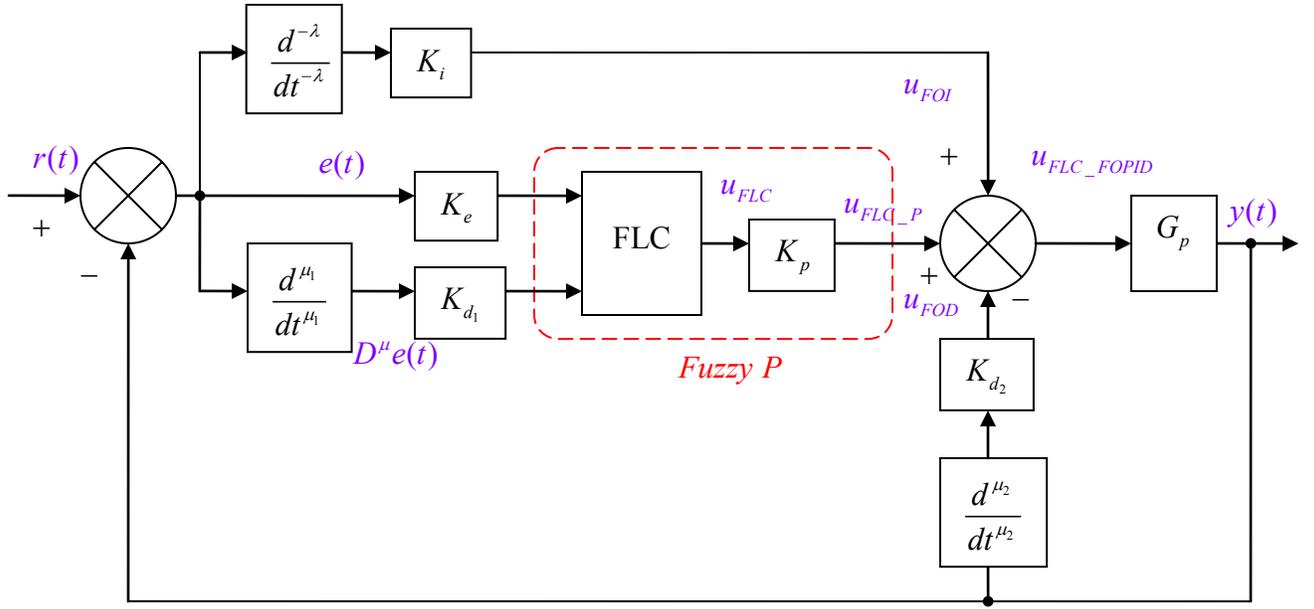

Fig. 3. Schematic diagram of FO fuzzy P+ID controller.

## 2.4. Fractional order fuzzy PI+D controller

Fuzzy logic based PI+D controller has been studied in [25]-[26]. Er and Sun [25] proposed GA based optimal tuning of PI+D controllers which has been extended in the present work with a fractional order enhancement of such controllers. Here, the derivative action in the feedback path not only gives smooth control action for sudden jump in set-point but also finely modulates the level of required compensation with additional flexibility of FO differentiator with order $\mu_2$.

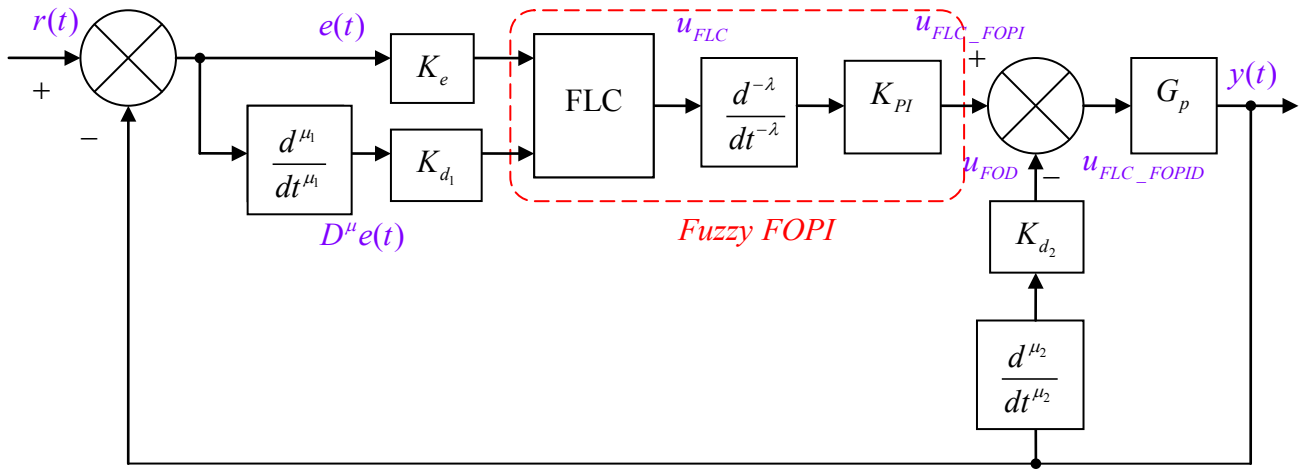

Fig. 4. Schematic diagram of FO fuzzy PI+D controller.

The control law of such a controller is given as

$$u_{FLC\_FOPID}(t) = u_{FLC\_FOPI}(t) - u_{FOD}(t)$$
$$= K_{PI} \cdot \frac{d^{-\lambda} u_{FLC}(t)}{dt^{-\lambda}} - K_{d_2} \cdot \frac{d^{\mu_2} y(t)}{dt^{\mu_2}} \tag{6}$$

### 2.5. Fractional order fuzzy PD+I controller

This structure is a parallel combination of FO fuzzy PD and FOPI controller. Since the present rule base uses two inputs, the single input FLC has not been introduced in the integral action unlike [27]. This fuzzy PD+I structure is just extension of the parallel or non-interacting structure of PID controllers in FO domain with the proportional and derivative actions being coupled and fuzzified. The control law of such fuzzy PD+I controller is given as

$$u_{FLC\_FOPID}(t) = u_{FLC\_FOPD}(t) + u_{FOI}(t)$$
$$= K_{PD} u_{FLC}(t) + K_i \cdot \frac{d^{-\lambda} e(t)}{dt^{-\lambda}} \tag{7}$$

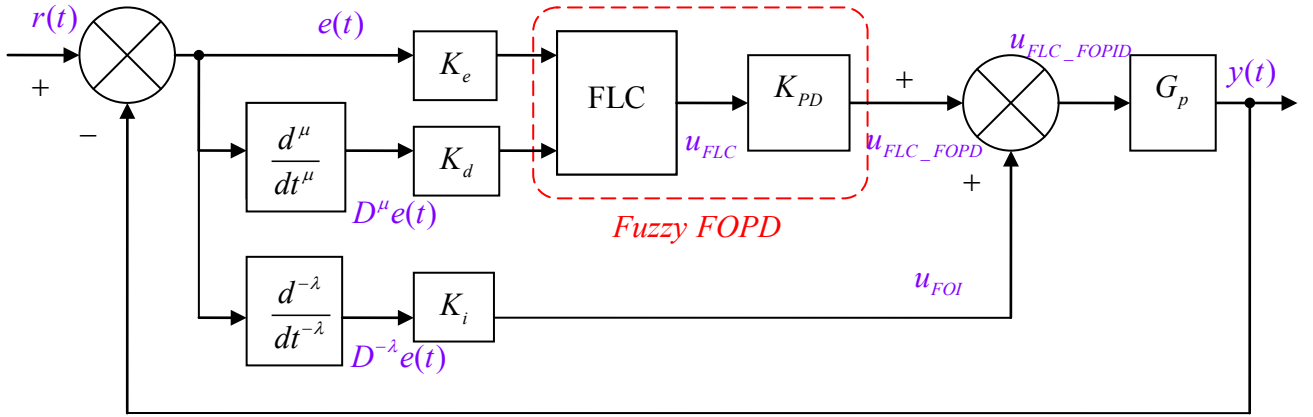

Fig. 5. Schematic diagram of FO fuzzy PD+I controller.

### 2.6. Details of the fuzzy inference and fractional order controller implementation

The fuzzy inference is the method by which the nonlinear mapping between the input and the output variables is established with the help of fuzzy logic. The process of fuzzy inferencing mainly comprises of—
  a) Fuzzy rule base
  b) Membership functions used in the rules
  c) Reasoning mechanism by the use of fuzzy logic operators
  d) Fuzzification and defuzzification opearations

The basic fuzzy logic based controller used here is based on a two dimensional rule base (Fig. 6) and triangular membership functions (Fig. 7), with 50% overlap. Other shapes of membership functions like Gaussian, Trapezoidal etc. could be used but it has been shown in Woo *et al.* [5] that the scaling factors affect the performance of the controller to a much greater extent than the shape of the membership functions. Also due to the fact that triangular membership function is the simplest one and easily implementable in hardware, it has been considered here for simulation studies. The set of rule-base has been divided into five groups as in Fig. 6 and its working principle has been illustrated in a detailed manner in Das *et al.* [28] in the context of FO fuzzy PID controller.

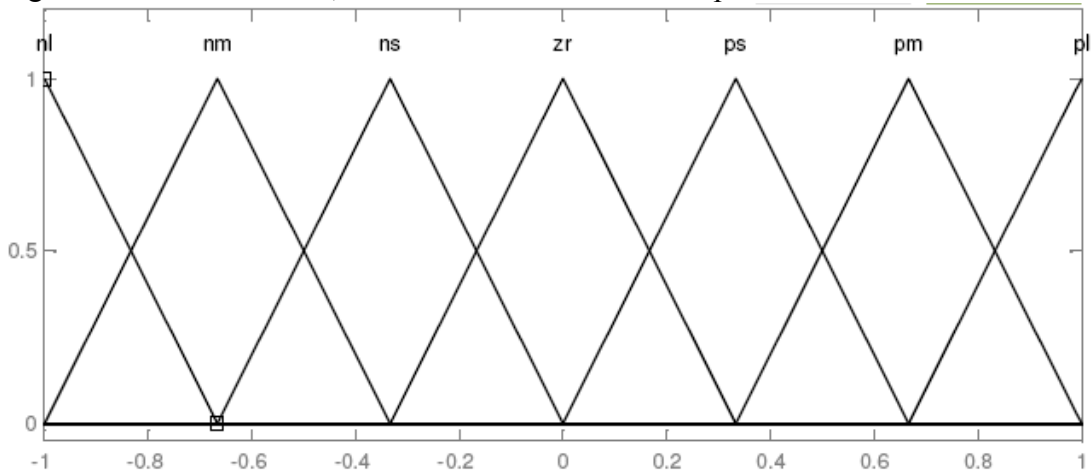

Fig. 6. Rule base for error, error derivative and FLC output.

Fig. 7. Membership functions for error, error derivative and FLC output.

Here, a Mamdani type inferencing is used with min type operator for implication and max type operator for rule aggregation. The error and its fractional derivative ($0 < \mu < 1$) is assumed to follow the rule base depicted in Fig. 6 composed of 49 ($7 \times 7$) rules. The acronyms NL, NM, NS, ZR, PS, PM and PL refer to Negative Large, Negative Medium, Negative Small, Zero, Positive Small, Positive Medium and Positive Large respectively. Since the linguistic variables dictate the granularity of the control action, more number of them could be used for better control resolution. But in these cases the rule base increases in the order of $n^2$ (where $n$ is the number of linguistic variables) and hence would be difficult to implement in real time hardware. The rule base is derived from [5] by incorporating expert knowledge in its design. However the fractional derivative of error signal is taken as the input here to the FLC instead of the integer order derivative which gives more flexibility while designing the controller as shown by Podlubny [13] for the case of $PI^\lambda D^\mu$ controller. The FLC outputs ($u_{FLC}$) in Fig. 1-5 is derived with the centre of gravity method for defuzzification. The control surface, describing the input-output relationship of the FLC is shown in Fig. 8. This essentially reflects the nonlinear mapping between the inputs (error and fractional rate of error) and the FLC output.

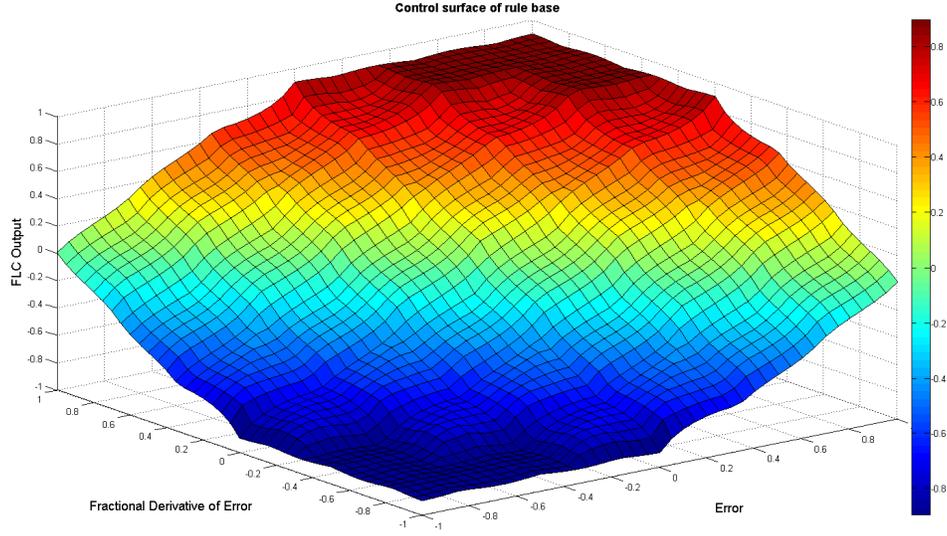

Fig. 8. Control surface for rule base.

Here, each of the proposed five classes of fractional order fuzzy PID families is designed with the same structure of fuzzy logic. It is arguable that different membership function type, number of membership functions or rule base could have changed the optimum performance of the corresponding closed loop control systems. But here we restricted with seven triangular membership functions and 49 rules for each of the fuzzy control structures to maintain a common platform for fair comparison of different configurations of input-output scaling factor and integro-differential orders using single and multi-objective optimization based approach.

It is clear that the each of the control laws for five different hybrid structures contains fractional order differentiation or integration of the error or FLC output respectively. The concept of fractional order PID type controller comes from the concept of fractional order differentiation and integration [1]-[2]. There are some popular definitions of fractional derivative like the Riemann-Liouville and Grunwald-Letnikov definitions. But, in the fractional order systems and control related literatures mostly the Caputo's fractional differentiation formula is referred. This typical definition of fractional derivative is generally used to derive fractional order transfer function models from fractional order ordinary differential equations with zero initial conditions. According to Caputo's definition the $\beta^{th}$ order derivative of a function $f(t)$ with respect to time is given by (8) and its Laplace transform can be represented as (9).

$$D^\beta f(t) = \frac{1}{\Gamma(m-\beta)} \int_0^t \frac{D^m f(t)}{(t-\tau)^{\beta+1-m}} d\tau, \quad \beta \in \mathbb{R}^+, m \in \mathbb{Z}^+ \tag{8}$$
$$m-1 \leq \beta < m$$

$$\int_0^\infty e^{-st} D^\beta f(t) dt = s^\beta F(s) - \sum_{k=0}^{m-1} s^{\beta-k-1} D^k f(0) \tag{9}$$

where, $\Gamma(\beta) = \int_0^t e^{-t} t^{\beta-1} dt$ is the Gamma function and $F(s) := \int_0^\infty e^{-st} f(t) dt$ is the Laplace transform of $f(t)$.

For control system analysis and design, it is often considered that the initial conditions of FO differential equations are zero to find out the transfer function representation of the linear FO dynamical system. With such an assumption the time domain operator $D^\beta$ can simply be represented in frequency domain as $s^\beta$. In this context, a negative sign in the derivative order $(-\beta)$ essentially implies a fractional integration operator. The classical FOPID or $PI^\lambda D^\mu$ controller [13] is therefore a weighted sum of such operators with extra degrees of freedom for tuning the weights (controller gains) along with the integro-differential order of the operators. This typical controller structure has five independent tuning knobs i.e. the three controller gains $\{K_p, K_i, K_d\}$ and two fractional order operators $\{\lambda, \mu\}$. For $\lambda = 1$ and $\mu = 1$ the FOPID controller structure reduces to the classical PID controller in parallel structure. In order to implement a fractional order control law, Oustaloup's band-limited frequency domain rational approximation technique is used in the present paper and also in most of the recent FO control literatures [1], [3]. In fact, the fractional control law with FO differ-integration can also be implemented using the Grunwald-Letnikov definition which is basically a finite difference approximation of fractional derivative with long memory behavior [1], [2]. But the rationale behind the choice of frequency domain rational approximation of FOPID controller is that it can be easily implemented in real hardware using higher order Infinite Impulse Response (IIR) type analog or digital filters, corresponding to each fractional order differ-integration in the FOPID controller.

On the other hand, the infinite dimensional nature of the fractional order differentiator and integrator in the FOPID controller structure creates hardware implementation issues in industrial application of FOPID controllers. However, few recent research results show that band-limited implementation of FOPID controllers using higher order rational transfer function approximation of the integro-differential operators gives satisfactory performance in industrial automation [29]. The Oustaloup's recursive approximation, which has been used to implement the integro-differential operators in frequency domain is given by the following expression [1], representing a higher order analog filter.

$$s^\beta \simeq K \prod_{k=-N}^{N} \frac{s + \omega'_k}{s + \omega_k} \tag{10}$$

where, the poles, zeros, and gain of the filter can be recursively evaluated as:

$$\omega_k = \omega_b \left(\frac{\omega_h}{\omega_b}\right)^{\frac{k+N+\frac{1}{2}(1+\beta)}{2N+1}}, \omega'_k = \omega_b \left(\frac{\omega_h}{\omega_b}\right)^{\frac{k+N+\frac{1}{2}(1-\beta)}{2N+1}}, K = \omega_h^\beta \tag{11}$$

Thus, any signal $f(t)$ can be passed through the filter (10) and the output of the filter can be regarded as an approximation to the fractionally differentiated or integrated signal $D^\beta f(t)$. In (10)-(11), $\beta$ is the order of the differ-integration, $(2N+1)$ is the order of the filter and $(\omega_b, \omega_h)$ is the expected fitting range. The advantage and disadvantages of Oustaloup's recursive filter based FOPID type controller simulation has been shown by Bayat with recommendation for improved performance given in [30]. Though for simplicity in the present study, 5$^{th}$ order Oustaloup's recursive approximation has been adopted to represent the integro-differential operators within a frequency band of $\omega \in \{10^{-2}, 10^2\}$ rad/sec.

Even with the truncation of infinite dimensional natures of FO operators with high order IIR filters, the obtained FOPID controllers are found to outperform classical PID structure in recent applications in process control especially networked control systems [31]. Thus there is always a trade-off between the complexity of the realization of the FOPID controller and the achievable accuracy. Further applications of fractional order PID type controllers in industrial automation related problems can be viewed in [31]-[37]. The focus of the present paper is to study the comparative efficacy of different hybrid structures of the nominal fractional order fuzzy PID controller, first introduced in [28] to handle oscillatory fractional order processes with dead time.

## 3. Time domain performance index based optimization for tuning of the family of fractional order fuzzy PID controllers
### 3.1. Formulation of the objective function

In the previous section several FLC based FOPID structures with their different control laws have been proposed. Now these structures need to be optimized to tune the control laws while meeting few control objectives incorporated as time domain performance indices. In the present study, the integral performance index ($J$) to be minimized by a suitable optimization algorithm has been taken as the weighted sum of ISTSE (Integral of Squared Time multiplied Squared Error) and ISDCO (Integral of Squared Deviation of the Controller Output) as follows:

$$J = \int_0^\infty \left[ w_1 \cdot t^2 e^2(t) + w_2 \cdot (u(t) - u_{ss})^2 \right] dt = (w_1 \times ISTSE) + (w_2 \times ISDCO) \quad (12)$$

Optimization result with objective function (12) produces the optimally tuned controller parameters (gains and integro-differential orders) in terms of low error index and control signal. The inclusion of the squared error term in the ISTSE penalizes the peak overshoot to a large extent. Also, the squared time multiplication term penalizes the error signal more at the later stages than at the beginning and hence results in a faster settling time. The squared deviation of the controller output is also included in $J$ so that the control signal does not become too large and result in actuator saturation and integral windup. In most fuzzy control problems at steady state the loop error and its rate goes to zero and as a result the steady state value of control signal ($u_{ss}$) also becomes zero which reduces the second part of expression (12) as simple squared control signal. The weights $w_1$ and $w_2$ have been incorporated in the objective function (12) to keep a provision for balancing the impact of the error and the control signal. In this case we have considered equal weights for the two objectives to be met by the controller i.e. minimal variation for controlled variable and manipulated variable as well. This implies that the minimization of the error index and the control signal are equally important. In the present study, a more stringent performance index (ISTSE) has been used unlike Integral of Time multiplied Absolute Error (ITAE) or Integral of Time multiplied Squared Error (ITSE) [28] to handle highly sluggish and oscillatory processes, though it might increase the control signal or create violent perturbation of the manipulated variable which has been taken into consideration during the optimization in (12). The impact of choosing different time domain performance index on fractional order optimum controller design has been shown in [28], [34].

The performance index (12) has two parts i.e. the error index and the integral of squared deviation of control signal which are then summed with some a priori decided weights for single objective optimization. These two parts may not be numerically of the same order for a particular case of controller and plant. Therefore the tracking and required controller effort will definitely vary for different structures though both of them have been equally weighted and their summation as a custom performance index has been minimized with the single objective optimization using real coded genetic algorithm. Also, in order to avoid a priori choice of the weights for the two parts of the objective function, optimal design trade-off between the two control objectives can be shown using multi-objective optimization which indicates that the designer has to sacrifice in one objective while improving the other. For the multi-objective NSGA-II based design trade-off comparison among different FO fuzzy controller structures the two set of conflicting objective functions are taken as:

$$\left.\begin{array}{l} J_1 = \int_0^\infty t^2 e^2(t) dt = ISTSE_{set-point} \\ J_2 = \int_0^\infty \left(u(t) - u_{ss}\right)^2 dt = ISDCO_{set-point} \end{array}\right\} \quad (13)$$

and

$$\left.\begin{array}{l} J_1 = \int_0^\infty t^2 e_{sp}^2(t) dt = ISTSE_{set-point} \\ J_3 = \int_0^\infty t^2 e_{ld}^2(t) dt = ISTSE_{load-disturbance} \end{array}\right\} \quad (14)$$

Here, $e_{sp}$ and $e_{ld}$ denotes the error signal for a particular choice of FO fuzzy controller structure and oscillatory fractional order process with unit change in the set-point and load-disturbance.

### 3.2. Single and Multi-objective optimization algorithm used for offline controller tuning

It has been shown by Pan and Das [3] that intelligent optimization based methods can be effectively used for such optimal fractional order controller design task. This has motivated to use a real coded genetic algorithm, which is a standard, widely used global optimization technique, for finding the optimal set of input-output SFs and integro-differential orders for the proposed family of FO fuzzy PID controllers as presented in section 2. GA is a stochastic optimization process inspired from Darwin's theory of evolution. It can be used to minimize a suitable objective function. Unlike other gradient based optimization techniques, the mechanism of GA is not based on the derivative of the objective function. Hence GA is suitable for searching global minima in highly nonlinear, rough and discontinuous functions and is less susceptible to be trapped in local minima like the gradient based methods. In GA a solution vector is initially randomly chosen from the solution search space and undergoes reproduction, crossover and mutation in each generation to produce a better population of solution vectors in the next generation. Reproduction implies that solution vectors with higher fitness values can produce more copies of themselves in the next generation. Crossover refers to information exchange based on probabilistic decisions between solution vectors. In mutation a small randomly selected part of a solution vector is occasionally altered, with a very small probability. This way the solution is refined iteratively until the objective function is minimized below a certain tolerance level or the maximum number of iterations are exceeded. The number of population members in GA is chosen to be 20. Another parameter called the elite count is used which reflects the number of fittest individuals in the current generation that would definitely be carried over to the next generation. This number is generally kept as a small percentage of the overall population to minimize the effect of dominance of initially obtained fitter individuals. In the present simulation the elite count is 2. Other than the elite genes, the rest of the population evolves through crossover and mutation which is dictated by the Crossover Ratio ($Cr$) and Mutation Ratio ($Mr$). The choice of $Cr$ and $Mr$ are problem dependent and must be chosen judiciously. In the present simulation $Cr=0.8$ and $Mr=0.2$ have been chosen which has proven to give good results for a wide variety of optimization problems [38]. Though similar optimization based FO controller designs have been attempted using other global optimization algorithms also, like Differential Evolution (DE) [31], Particle Swarm Optimization (PSO) [20], [3] etc. the present paper uses GA, as a popular method of optimization based controller finding [39], [28], [3]. Also, the GA has been run multiple times to ensure that the true global minima has been found in the optimization process and the best results in terms of the lowest cost function have been reported along with the optimum decision variables, in the sub-sequent sections showing simulation examples with the hybrid fuzzy FOPID controllers.

Also, each of the proposed FO fuzzy PID controller structure has a performance limit. For example, it is desirable that the controller results in a very fast settling time and also simultaneously produces lower control signal. But these are contradictory objectives [40]-[41] and design considerations to minimize one would definitely make the other larger. Thus an effective way of comparing the limits of the various controller structures is to study the design trade-offs between the different contradictory criteria like set point tracking vs. disturbance rejection capability and fast error minimization vs. required control signal etc. For this study a multi-objective evolutionary algorithm known as NSGA-II [42] has been employed and the non-dominated Pareto front is obtained for the each of the controller structures. The NSGA-II algorithm operates in a similar fashion as the single objective GA employing crossover, mutation and reproduction. Additionally, since the solutions give multiple objective values, they are selected based on their fitness as well as spread on the non-dominated front [42]. A non-domination rank and crowding distance is assigned to the solutions and they are sorted using a fast non-dominated sorting algorithm [42]. This helps in finding solutions very close to the true optimal Pareto frontier and helps in preserving the diversity of the solutions to obtain the whole length of the Pareto front. In the present simulation, the number of population is increased to 100 as compared to the single objective GA. The crossover and mutation factors and elite count are kept as in the previous single objective GA cases and the Pareto fraction (denoting the fraction of total solutions that should be on the Pareto front) as 0.7.

### *3.3. Constraints imposed to guarantee closed loop stability and faster convergence*

Analytical stability issues of such fuzzy PID controllers have been studied in [18]. Also, with fractional order models in the loop the scope of classical Matignon's stability [43], [1] or Kharitonov like FO robust stability checking criterion [44] can't be applied directly due to the inclusion of high nonlinearity of the FLC in the control loop. As an alternate to this an equivalent criteria has been imposed within the optimization algorithm to guarantee the closed loop stability of such controllers having fuzzy modules [20] and FO elements [28]. It is obvious that the analytical stability criteria in the optimization algorithm can be alternately incorporated by removing the unstable modes within the search space that can be easily found by observing the magnitude of the integral performance index (12) itself for an intermediate guess solution vector of GA and NSGA-II. The optimization algorithm is formulated such that it automatically rejects infeasible values of the input-output SFs that might cause instability of the closed loop system. A heavy penalty is imposed in the objective function when the ISTSE value becomes large for some choice of solution variables in the unstable region like that adopted in [28]. Thus the algorithm automatically steers off these areas in the search space and converges to areas which give better values for error index and the control signal also.

In the genetic algorithm, constraints are also imposed on the lower and upper bound of the solution variables. All the input scaling factors for the FLC variants are constrained to be within $[0,1]$ to fix the universe of discourse for fuzzy inference. The output scaling factor is constrained also to vary in the interval $[0,40]$ since they are equivalent to the controller gains to avoid the difficulty in practical implementation issues for high controller gains. All the integro-differential orders are constrained to vary between $[0,2]$, beyond which the closed loop system generally shows instability. It has been found that such constrained version of single or multi-objective GA gives better performance compared to its unconstrained version [38], in terms of faster convergence, quality of obtained solution and avoidance of local minima.

### 4. Simulation and Results:

The main focus of the paper is to show performance comparison of different fractional order fuzzy PID controller structures with respect to three different control objectives for fractional order plants showing oscillatory open loop response. Recently it has been shown in Das *et al.* [34] that many higher order systems can be compactly and accurately represented by fractional order models known as Non-Integer Order Plus Time Delay (NIOPTD). Since, it is well known that fractional order controllers are the best means of controlling such systems [1], such single and multi-objective optimization based comparison may lead to meaningful conclusions over the state-of-art techniques.

In this section we tested the comparative merits of the proposed hybrid FO fuzzy PID controllers handling few complicated processes. The optimal controller gains and integro-differential orders are found out by the GA based optimization presented in the previous section. FO controller design for First Order Plus Time Delay (FOPTD) integer order systems has been extensively studied in Chen *et al.* [35] and Saha *et al.* [36], among many other literatures. In this paper, three fractional order systems have been tuned with the family of FO fuzzy controllers which shows highly oscillatory open loop responses but has different values of normalized dead-time ratio [35], [39]. The generalized template representing the three test-processes is given as (15) known as NIOPTD-I [34].

$$G_p(s) = \frac{K}{Ts^\alpha + 1} e^{-Ls} \quad (15)$$

A wide variety of oscillatory open loop processes with higher order dynamics can be effectively modeled using the template in (15) which motivates us to study control of this typical structure. It has been seen that such a fractional order system shows highly sluggish and oscillatory open loop time response for $(0 < \alpha < 1)$ and $(1 < \alpha < 2)$ respectively. The processes chosen for optimum hybrid fuzzy PID controller tuning are similar to the FOPTD processes in [36] with different time delay to time constant ratio ($L/T$) often found in process control literatures [39] and additionally the order of the fractional order systems $(\alpha = 1.5)$ producing oscillatory open loop response.

### *4.1. Control of lag dominated (L<<T) oscillatory fractional order (1<α<2) process*

Let us consider, the lag dominated process in [36] with oscillatory open loop dynamics

$$G_{p1}(s) = \frac{1}{1.11s^{1.5} + 1} e^{-0.105s} \quad (16)$$

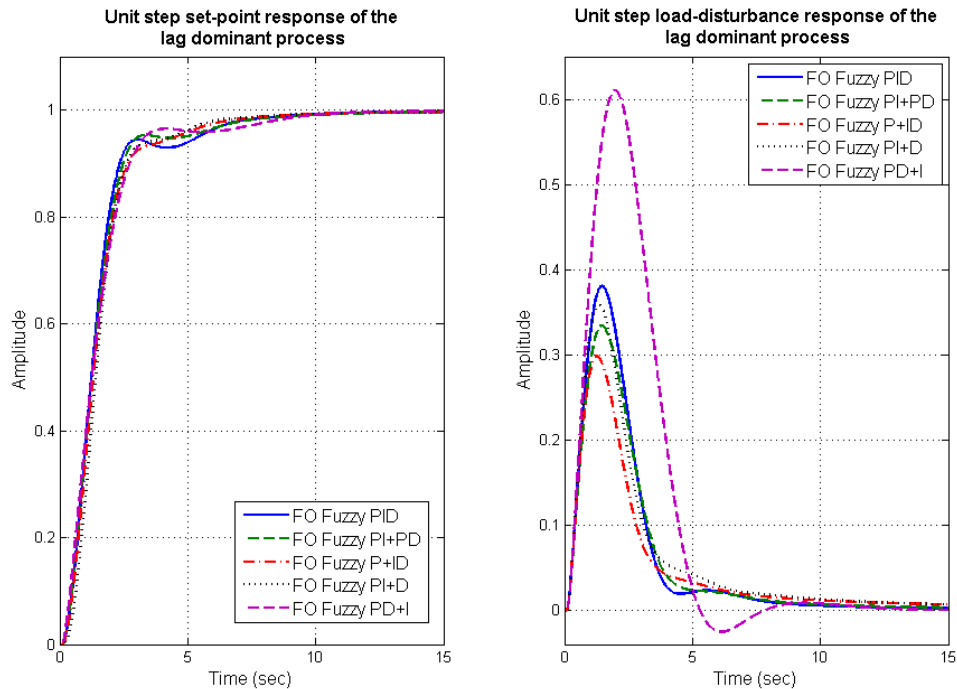

Fig. 9. Response of the lag dominant process with unit step change in set-point and load disturbance.

The GA based optimum tuning results of the process (16) with the proposed family of FO hybrid fuzzy PID controllers have been reported in Tables 1-5 and the corresponding closed loop performances are shown in Fig. 9-10. In Fig. 9, it is clear that the unit step set-point response is almost closer to each other for all the different controller structures. However, the load disturbance rejection of the FO fuzzy P+ID structure is better for the lag dominated process than the other ones. The FO fuzzy PD+I controller has the worst load disturbance rejection amongst the rest. In Fig. 10, the control signal for unit set point change is nearly similar for all the other controllers, except the FO fuzzy PI+D controller which has a higher output at the onset of the step input.

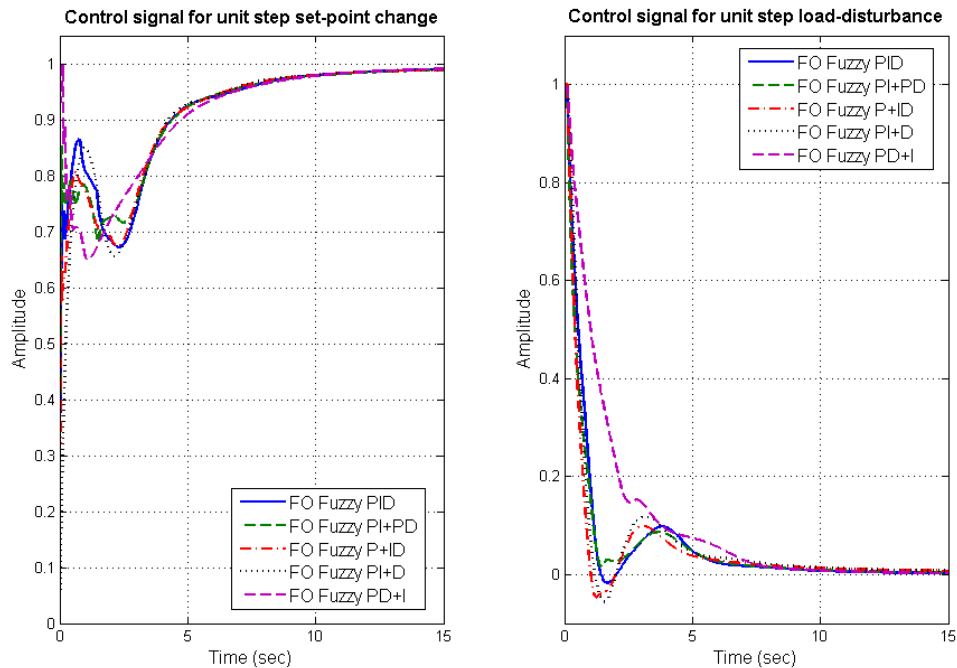

Fig. 10. Control signal for the lag dominant process with unit step change in set-point and load disturbance.

Table 1
Optimal tuning result for FO fuzzy PID controller

| Process | $J_{min}$ | Optimal FO fuzzy PID Controller Parameters | | | | | |
|---|---|---|---|---|---|---|---|
| | | $K_e$ | $K_d$ | $K_{PI}$ | $K_{PD}$ | $\lambda$ | $\mu$ |
| $G_{p1}$ | 38.20247 | 0.887976 | 0.63353 | 1.417276 | 0.820367 | 0.959188 | 0.994714 |
| $G_{p2}$ | 7.630405 | 0.098897 | 0.102872 | 0.728721 | 0.787448 | 0.998849 | 0.992102 |
| $G_{p3}$ | 39.6631 | 0.666385 | 0.214853 | 0.801473 | 0.321055 | 0.998524 | 0.288179 |

Table 2
Optimal tuning result for FO fuzzy PI+PD controller

| Process | $J_{min}$ | Optimal FO fuzzy PI+PD Controller Parameters | | | | | | | |
|---|---|---|---|---|---|---|---|---|---|
| | | $K_{e_1}$ | $K_{d_1}$ | $K_{PI}$ | $K_{e_2}$ | $K_{d_2}$ | $K_{PD}$ | $\lambda$ | $\mu$ |
| $G_{p1}$ | 38.17563 | 0.957059 | 0.74568 | 1.506117 | 0.725838 | 0.872039 | 0.882793 | 0.932188 | 0.982342 |
| $G_{p2}$ | 3.752172 | 0.177834 | 0.016532 | 0.636613 | 0.299998 | 0.765192 | 0.287097 | 0.976782 | 0.810926 |
| $G_{p3}$ | 39.64602 | 0.848295 | 0.209849 | 0.843522 | 0.295589 | 0.209216 | 0.487242 | 0.971632 | 0.436048 |

Table 3
Optimal tuning result for FO fuzzy P+ID controller

| Process | $J_{min}$ | Optimal FO fuzzy P+ID Controller Parameters | | | | | | | |
|---|---|---|---|---|---|---|---|---|---|
| | | $K_e$ | $K_{d_1}$ | $K_p$ | $K_{d_2}$ | $K_i$ | $\lambda$ | $\mu_1$ | $\mu_2$ |
| $G_{p1}$ | 38.1687 | 0.339126 | 0.81547 | 0.594271 | 1.924765 | 1.806937 | 0.882179 | 0.973166 | 0.177353 |
| $G_{p2}$ | 3.631472 | 0.007836 | 0.288275 | 0.650441 | 0.131799 | 0.17253 | 0.973567 | 0.769968 | 0.05902 |
| $G_{p3}$ | 39.69599 | 0.64044 | 0.094509 | 0.301722 | 0.161946 | 0.657659 | 0.972741 | 0.998061 | 0.00964 |

Table 4
Optimal tuning result for FO fuzzy PI+D controller

| Process | $J_{min}$ | Optimal FO fuzzy PI+D Controller Parameters | | | | | | |
|---|---|---|---|---|---|---|---|---|
| | | $K_e$ | $K_{d_1}$ | $K_{PI}$ | $K_{d_2}$ | $\lambda$ | $\mu_1$ | $\mu_2$ |
| $G_{p1}$ | 38.21658 | 0.658696 | 0.328859 | 2.02627 | 1.314265 | 0.883782 | 0.707495 | 0.432665 |
| $G_{p2}$ | 6.67324 | 0.435695 | 0.240776 | 0.379578 | 0.314335 | 0.873519 | 0.59048 | 0.753619 |
| $G_{p3}$ | 39.89151 | 0.712596 | 0.20361 | 1.06411 | 0.220181 | 0.940606 | 0.607729 | 0.429407 |

Table 5
Optimal tuning result for FO fuzzy PD+I controller

| Process | $J_{min}$ | Optimal FO fuzzy PD+I Controller Parameters | | | | | |
|---|---|---|---|---|---|---|---|
| | | $K_e$ | $K_d$ | $K_i$ | $K_{PD}$ | $\lambda$ | $\mu$ |
| $G_{p1}$ | 38.22424 | 0.207274 | 0.59619 | 0.639649 | 1.039919 | 0.983022 | 0.599213 |
| $G_{p2}$ | 3.297377 | 0.056807 | 0.211725 | 0.113836 | 0.828508 | 0.989822 | 0.723279 |
| $G_{p3}$ | 39.67555 | 0.344379 | 0.5251 | 0.626799 | 0.33055 | 0.96105 | 0.28574 |

### *4.2. Control of balanced lag and delay type (L≈T) oscillatory fractional order (1<α <2) process*

The FO oscillatory process considered next has closer values of the time constant and delay having the following transfer function [36]:

$$G_{p2}(s) = \frac{5}{1.5s^{1.5}+1}e^{-s} \tag{17}$$

The corresponding optimum controller performances are shown in Fig. 11-12. As seen from the $J_{min}$ values in Tables 1-5, the FO fuzzy PI+PD, FO fuzzy P+ID and FO fuzzy PD+I controllers have good performances. The load disturbance response of the FO fuzzy PI+PD is the worst, in spite of having nice set-point tracking characteristics. The FO fuzzy P+ID controller has a smaller overshoot than the FO fuzzy PD+I controller but undershoot of the former is greater than the later in the load disturbance rejection plot. The control signal with the FO fuzzy PD+I structure has lesser magnitude of oscillation than FO fuzzy P+ID controller. The FO fuzzy PI+D has a large overshoot in the unit set-point change, but has a very good load disturbance response.

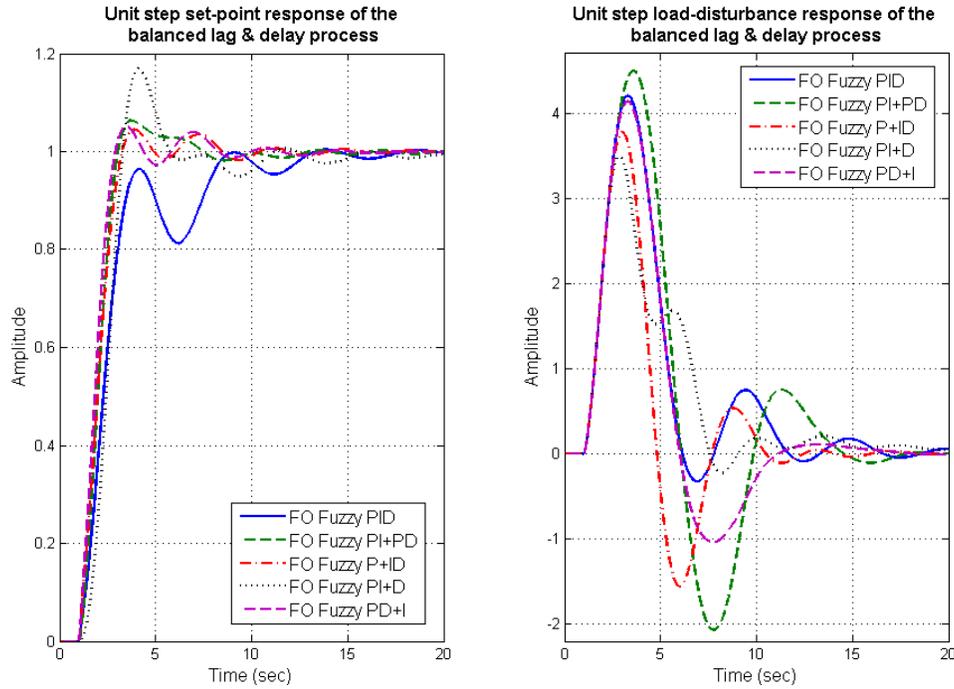

Fig. 11. Response of the balanced lag and delay process with unit step change in set-point and load disturbance.

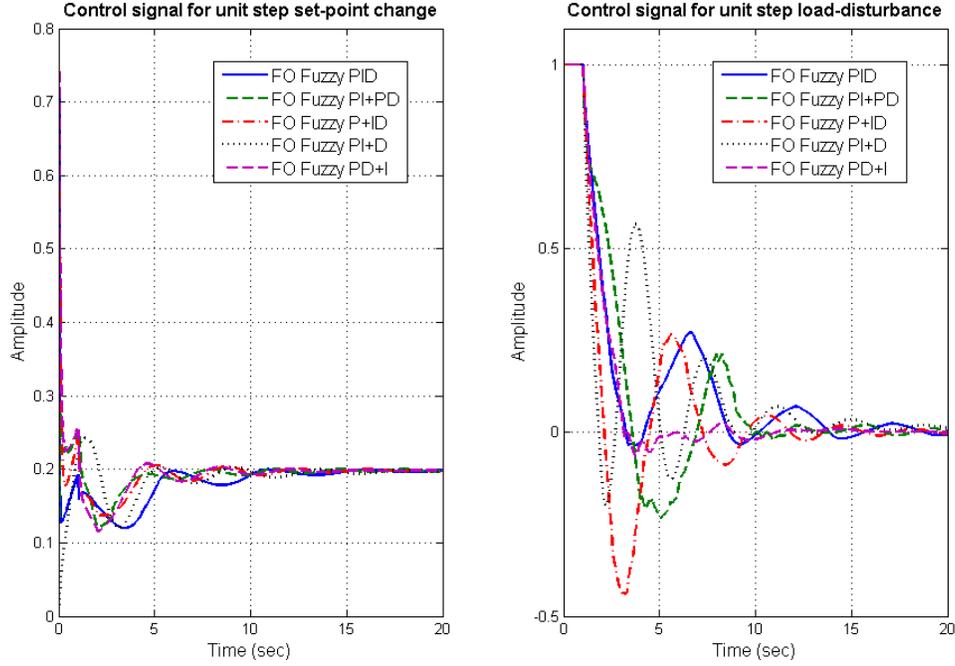

Fig. 12. Control signal for the balanced lag and delay process with unit step change in set-point and load disturbance.

### *4.3. Control of delay dominated (L>>T) oscillatory fractional order (1<α<2) process*

For the following delay dominant oscillatory FO process [36] the optimum fuzzy FOPID controllers are investigated next:

$$G_{p3}(s) = \frac{1}{0.05s^{1.5}+1} e^{-s} \quad (18)$$

It is clear from the process model that although the process has small dead-time (1 sec) but the relative dead time ($\tau$) is large having the following expression.

$$\tau = \frac{L}{L+T} \quad (19)$$

It has been shown in Chen *et al.* [35] that performance of FO controllers does not solely depend on the process time constant or delay but on the relative dead-time ($\tau$). The corresponding closed loop responses for process (18) and optimal FO hybrid fuzzy PID controllers are shown in Fig. 13-14. It is evident that for delay dominated processes the controllers have almost similar unit set-point response and load disturbance rejection response. However the control signal with unit step set-point and load disturbance is somewhat large and oscillatory respectively for the FO fuzzy PI+D controller.

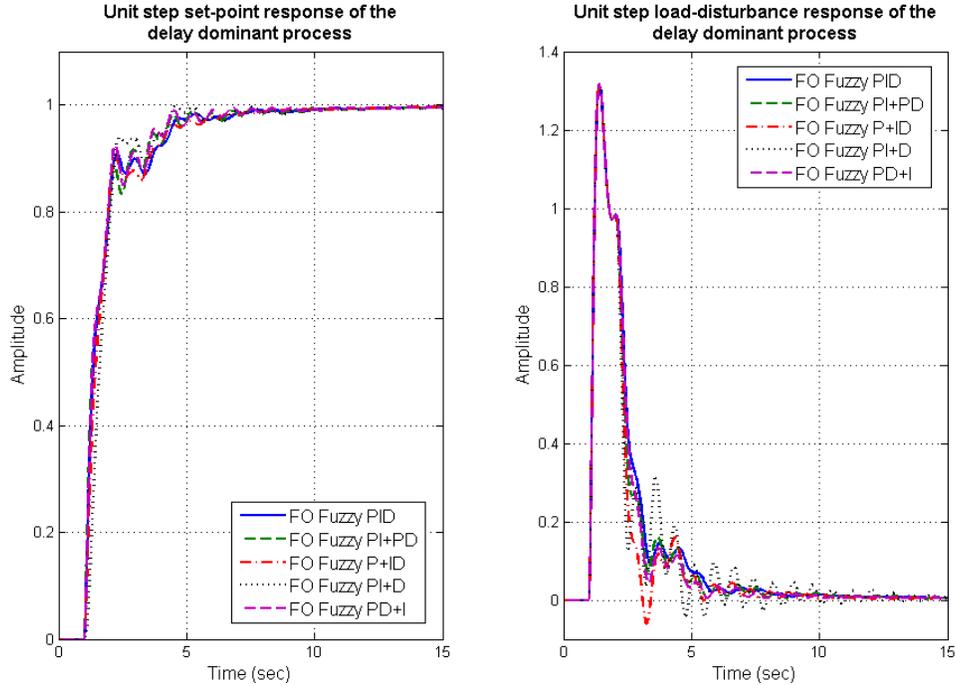

Fig. 13. Response of the delay dominant process with unit step change in set-point and load disturbance.

Thus, from the simulation results for three FO oscillatory test processes with different delay to lag ratio, it can be summarized that overall the FO fuzzy PI+D controller works well for the balanced lag and delay type oscillatory fractional order processes. However, the overall recommendation is higher for FO fuzzy P+ID structure to handle lag dominated oscillatory fractional order processes. It is obvious that for small magnitude of relative dead time the FO derivative action implemented on the process variable instead of the loop error produces better response as can also be found in the classical case [8]. For delay dominated oscillatory fractional order processes having larger value of $\tau$, these controllers does not perform well and the FO Fuzzy PID controller shows a superior performance for compensation of processes with larger relative dead-time.

Table 6
Recommended structure of the FO hybrid fuzzy controllers, found using single objective optimization, for handling different type of processes with different control objectives

| Type of Process | Best Controller Structure for Different Control Performance | | |
|---|---|---|---|
| | Set-point tracking | Load disturbance rejection | Small control signal |
| lag-dominant | almost similar for all structures | FO fuzzy P+ID | FO fuzzy PID |
| balanced lag and delay | FO fuzzy P+ID | FO fuzzy PI+D | FO fuzzy PD+I |
| delay dominant | almost similar for all structures | FO fuzzy PD+I | FO fuzzy PD+I |

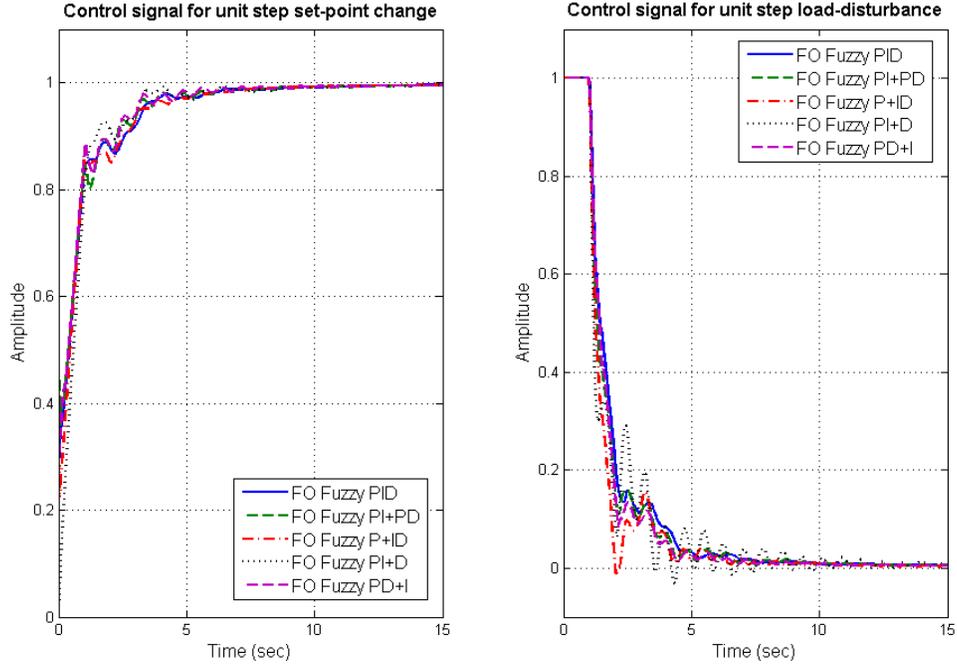

Fig. 14. Control signal for the delay dominant process with unit step change in set-point and load disturbance.

It is fact that the single objective optimization with (12) has been carried out only with set-point changes and additionally we have shown the load-disturbance characteristics of such optimum set-point tracking based FO fuzzy PID controllers. Since, modulating the maximum magnitude of sensitivity function to control the load disturbance characteristics like that in [35] are difficult and more mathematically involved for highly nonlinear controllers as in our case, we restricted our study on various performance comparison with the family of FO fuzzy PID controllers for optimum set-point based tuning only. The variation in control signal or manipulated variable has also been taken into consideration for unit set-point change only in the optimization based controller tuning process. The summary of the comparative performances of the family of fuzzy FOPID controllers are presented in Table 6 for three classes of oscillatory fractional order processes with various levels of relative dead-time. The proposed family of FO hybrid fuzzy PID controller structure is believed to dominate future process control industries over present day's fuzzy PID controllers, if the hardware implementation issues can be circumvented for both the fuzzy [45] and fractional differ-integral modules [1]-[2]. In addition to the recommended structure, it is interesting to see the achievable design trade-offs for different controller structure which requires multi-objective formulation of the controller tuning problem using different conflicting objectives (13) and (14).

*4.4. Multi-objective simulation comparison for three types of FO processes and five hybrid FO fuzzy PID controllers*

It is well known that a single controller structure cannot give good results for all design specifications. For specific applications, different controller structures would give a trade-off solution among conflicting design objectives. Hence for effective comparison of different controller structures it is essential to know the limits of performance of each of the individual controllers for contradictory design specifications. In [40]-[41], a similar approach has been taken to compare the efficacy of the fractional order PID controller vis-à-vis the integer order one. In the present case, two specific set of contradictory objectives are considered:
   a) Set point tracking vs. controller effort represented by (13)
   b) Set point tracking vs. load disturbance rejection performance represented by (14)

The reason for considering these two as contradictory objective functions can be briefly explained as follows. To achieve a faster set point tracking, it is essential that the controller gains should be higher. In other words, the controller must be able to exert much more control action on the process so that it settles in a short amount of time. However, the control signal should ideally be smaller to prevent actuator saturation and minimize the cost associated with sizing of a larger actuator. It can be inferred that both these objectives of small control signal, as well as faster set point tracking cannot be ideally obtained by a single objective optimization based controller. Thus, there would exist a range of values for the tuning parameters of the controller, where the controller would show good set point tracking at the cost of higher control signal and vice-versa.

A similar argument can be drawn for the conflicting objectives of set point tracking and load disturbance rejection as well. In [46], it has been explained that for any PID controller design, there are two objectives servo (command following) and regulatory (load disturbance rejection). Different tuning rules have been developed over the past few decades for any one of the criterion and it has been shown that these rules designed for one criterion might not give good results for the other one. Servo/regulatory design trade-off based PID tuning is illustrated in [47].

Fig. 15 shows the Pareto optimal trade-off between tracking and load disturbance rejection for oscillatory FO lag dominant process (16) with different controller structures. These solutions on the Pareto front represent only the non-dominated ones, i.e. in the NSGA-II run, other solution could be found which had a lower value of both the objectives together. Hence in a way, these represent the limits of each controller's performance. Among all the different proposed FO fuzzy controllers, it can be observed from Fig. 15 that the FO fuzzy P+ID controller gives the best performance with respect to the control signal and set point tracking. This is because, all the other Pareto fronts lie inside the concave region of this one, indicating that the FO fuzzy P+ID controller has lower values of the objective functions than the others. Also the spread of the Pareto for this controller is much more than the others. This indicates that the FO fuzzy P+ID controller is able to give more diverse solutions and the designer can choose one particular solution from the whole range, depending on his specific requirements. For example, if in the design problem the control cost is very expensive, then the solution which has a lower control cost should be chosen, although there would be a decrease in the performance of set point tracking. Fig. 16 shows the Pareto fronts between the tracking and load disturbance rejection for lag dominant processes with various FO fuzzy controller structures. In this case the FO fuzzy PD+I controller gives the best Pareto front.

Most of the other Pareto fronts are located in the concave region of this one, representing inferior solution sets.

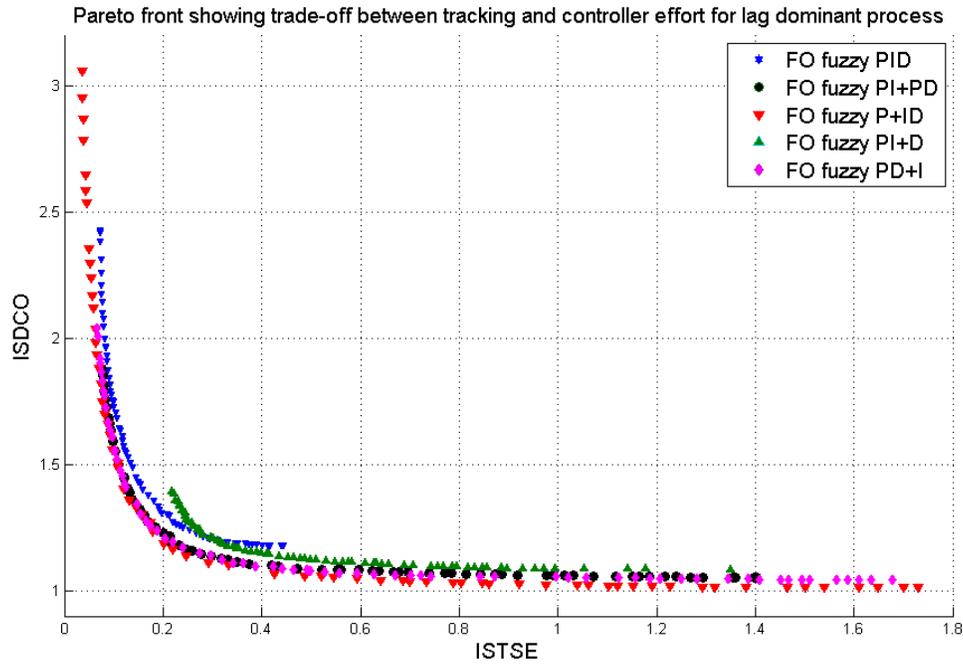

Fig. 15. Pareto front showing trade-off between tracking and controller effort for lag dominant process using various FO Fuzzy controller structures.

Fig. 17 and Fig. 18 show the Pareto fronts for the balanced lag and delay process (17). The former shows the tradeoff between the tracking and the controller effort, while the latter shows the tradeoff between the tracking and disturbance rejection. The two best controller structures as can be seen from Fig. 17 are the FO fuzzy P+ID and the FO fuzzy PD+I. However, the Pareto fronts of these two controller structures are intersecting in nature. Thus no clear winner emerges in this case. The FO fuzzy PD+I structure gives the best tracking performance at the cost of higher control signal. On the other hand the FO fuzzy P+ID gives a lower value of the control signal at the expense of sluggish tracking performance. In Fig. 18 the FO fuzzy PI+D gives the best load disturbance rejection properties with sufficiently good set point tracking. The FO fuzzy PI+PD controller gives better performance in set point tracking than this one, but the load disturbance rejection is significantly poor.

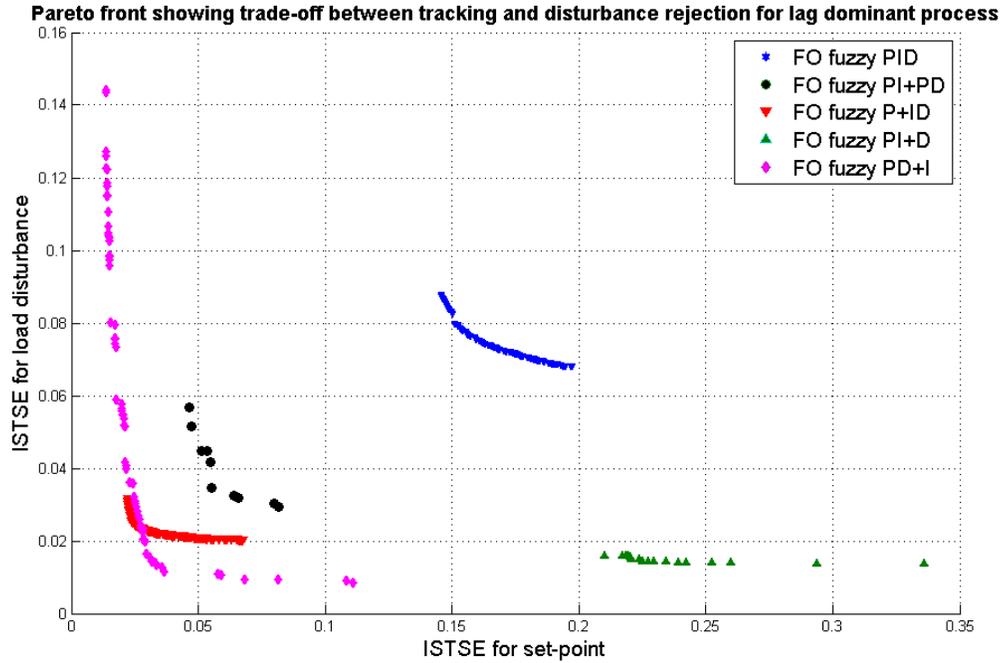

Fig. 16. Pareto front showing trade-off between tracking and disturbance rejection for lag dominant process using various FO Fuzzy controller structures.

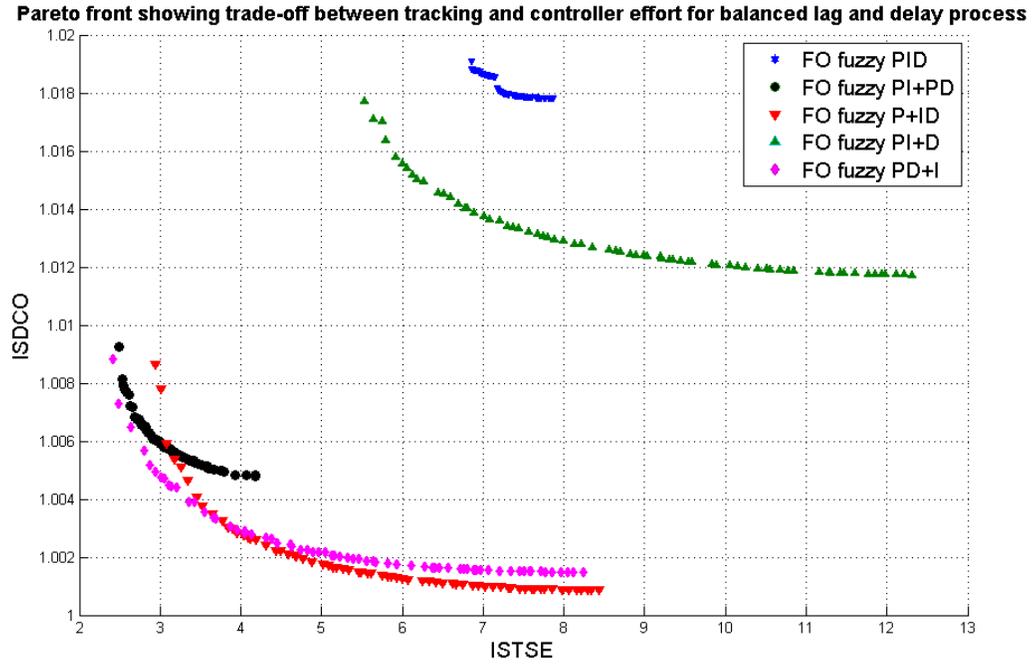

Fig. 17. Pareto front showing trade-off between tracking and controller effort for balanced lag and delay process using various FO Fuzzy controller structures.

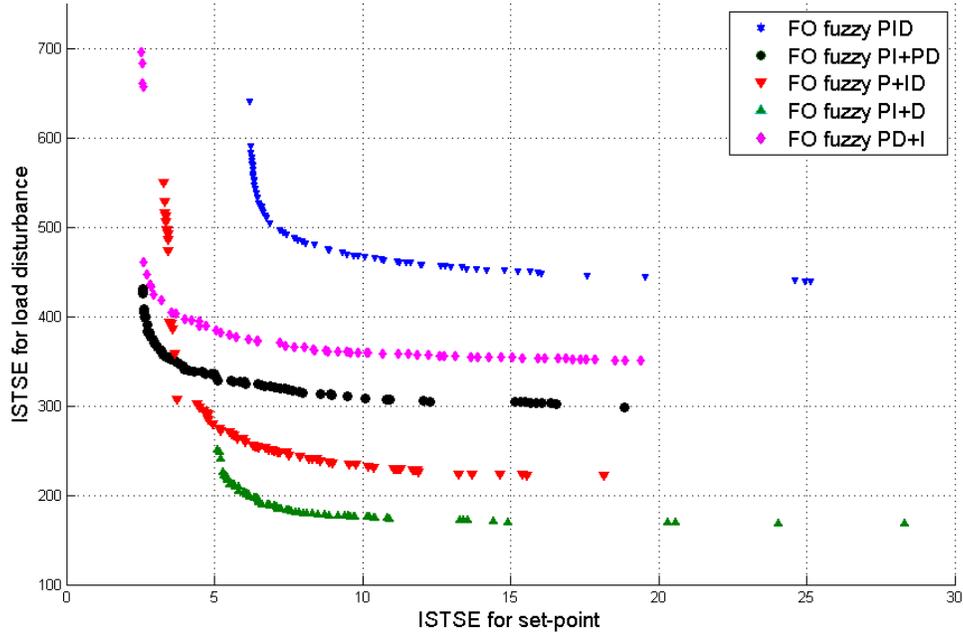

Fig. 18. Pareto front showing trade-off between tracking and disturbance rejection for balanced lag and delay process using various FO Fuzzy controller structures.

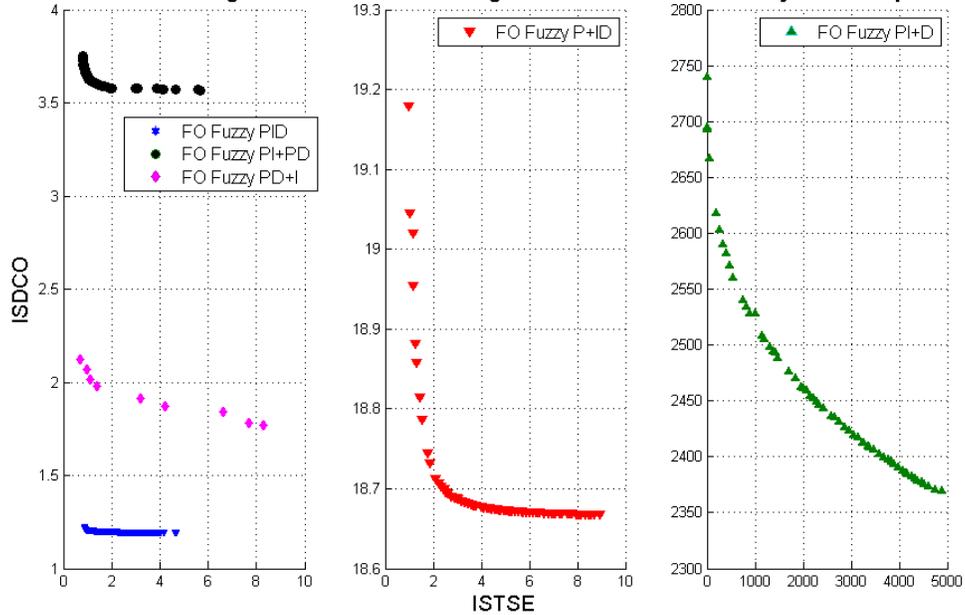

Fig. 19. Pareto fronts showing trade-off between tracking and controller effort for delay dominant process using various FO Fuzzy controller structures.

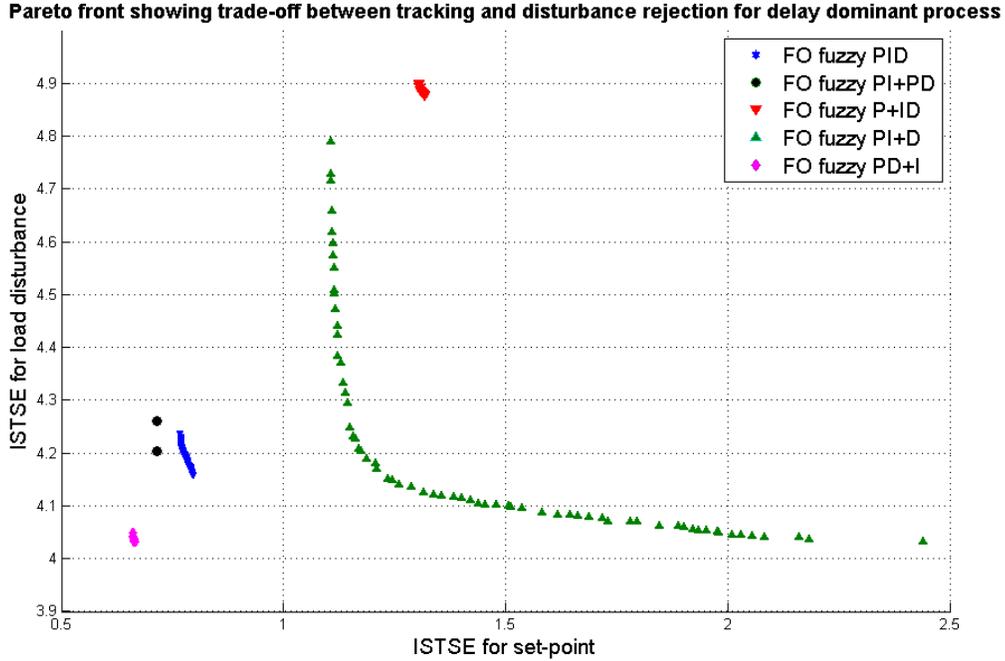

Fig. 20. Pareto fronts showing trade-off between tracking and disturbance rejection for delay dominant process using various FO Fuzzy controller structures.

Fig. 19 and Fig. 20 show the Pareto fronts for the delay dominant process (18). The former shows the tradeoff between the tracking and the controller effort, while the latter shows the tradeoff between the tracking and disturbance rejection. It can be seen that the controllers show wide variation in performance, for handling this kind of processes. From Fig. 19 it can be seen that the FO fuzzy PID gives the best results for lower control signal and set point tracking. The FO Fuzzy PI+D structure is the worst performer and is not suitable for processes with large time delay.

Table 7
Recommended structure of the FO hybrid fuzzy controllers, found using multi-objective optimization, for handling different type of processes with conflicting control objectives

| Type of Process | Best Controller Structure for Different Control Performance ||
| --- | --- | --- |
| | Set-point tracking vs. Small control signal | Set-point tracking vs. Load disturbance rejection |
| lag-dominant | FO fuzzy P+ID | FO fuzzy PD+I |
| balanced lag and delay | FO fuzzy P+ID and FO fuzzy PD+I | FO fuzzy PI+D |
| delay dominant | FO fuzzy PID | FO Fuzzy PD+I |

From Fig. 20 it can be seen that the FO fuzzy PD+I gives the best results for the conflicting objectives of set point tracking and load disturbance rejection. However the FO fuzzy PI+D controller has a more Pareto spread, indicating diverse solutions, although the performance is lesser than that of the FO fuzzy PD+I controller. It is possible to show the time domain responses for some representative solutions on the Pareto fronts, as has been done in other literatures like [40]-[41], to enunciate the effect of the tradeoff obtained amongst the different solutions. However, due to paucity of space, this illustration is not done in the present paper. The Pareto fronts are a sufficient indication of the limits achievable by each controller structure and can be used to understand the time domain characteristics of the plants tuned with each specific controller. The recommended controller structure based on the comparison of the design trade-offs between conflicting objectives are given in Table 7.

**5. Conclusion**

Several decomposed hybrid structures of the fractional order fuzzy PID controller have been proposed in this paper along with its optimal time domain tuning. Genetic algorithm based tuning of such controllers are attempted by minimizing a weighted summation of loop error index (with higher powers of time and error term to ensure fast and smooth tracking) and integral of squared deviation of control signal as the performance index. Three different classes of oscillatory fractional order processes (viz. lag dominant, balanced lag-delay and delay dominant) have been attempted to control with the proposed family of hybrid fractional order fuzzy PID controllers. A multi-objective optimization algorithm has been employed to study the performance trade-offs for different FO fuzzy PID family of controller structures. For different set of control objectives and different class of processes, various recommendations have been made and the best controller structures have been identified. Comparative merits and recommendations among the proposed family of FO fuzzy controllers are shown in terms of good set-point tracking, good load-disturbance rejection performance and minimal variation of the manipulated variable or control signal, for the compensation of processes with different levels of relative dead-time. The comparison of different FO fuzzy PID structures for different FO oscillatory processes shows how well a particular controller is capable of executing a specific task or two conflicting tasks which are generally incorporated in terms of various integral performance indices. Future scope of work can be directed towards finding analytical stability criteria for the family of proposed fuzzy FOPID controller structures.